\documentstyle[12pt]{article}

\newcommand{\sect}[1]{\setcounter{equation}{0}\section{#1}}
\newcommand{\subsect}[1]{\subsection{#1}}

\def\be{\begin{equation}}
\def\ee{\end{equation}}
\def\bea{\begin{eqnarray}}
\def\eea{\end{eqnarray}}

\def\jj{K} 
\def\pu{{H}} 
\def\pd{{P}}
\def\cu{{C_1}}
\def\cd{{C_2}}
\def\dd{D}

\def\>#1{{\bf #1}}                 
\def\1{\'{\i}}                           
\def\R{{\rm I\kern-.2em R}}

\begin{document}

\thispagestyle{empty}

\hfill \quad \ 
\ 
\vspace{2cm}

\begin{center}
{\LARGE{\bf{(2+1) null-plane quantum Poincar\'e group}}} 
  
{\LARGE{\bf{from a factorized universal $R$-matrix}}} 
  
\end{center}

\bigskip\bigskip\bigskip

\begin{center} Angel Ballesteros and  Francisco J. Herranz
\end{center}

\begin{center} {\it {  Departamento de F\1sica, Universidad
de Burgos} \\   Pza. Misael Ba\~nuelos  \\
E-09001, Burgos, Spain}
\end{center}

\bigskip\bigskip\bigskip

\begin{abstract}
The non-standard (Jordanian) quantum deformations of  $so(2,2)$ and
(2+1) Poincar\'e algebras are constructed by starting from a 
quantum $sl(2,\R)$ basis such that simple factorized 
expressions for their
corresponding universal $R$-matrices are obtained. As  an
application, the null-plane quantum (2+1) Poincar\'e
Poisson-Lie group is quantized by following the FRT
prescription. Matrix and differential representations of this
null-plane deformation are presented, and the influence of the
choice of the basis in the resultant $q$-Schr\"odinger equation
governing the deformed null plane evolution is commented.
\end{abstract}

\newpage


\sect {Introduction}

Among quantum deformations of Poincar\'e algebra we find
three remarkable Hopf structures. Two of them are   obtained in
a natural way within a purely kinematical framework encoded
within the usual Poincar\'e basis.
 They are the well-known $\kappa$-Poincar\'e algebra
\cite{Luk,Gill}  where the deformation parameter can  be
interpreted as a fundamental time scale and a $q$-Poincar\'e
algebra \cite{Afin} where the quantum parameter is a
fundamental length. On the other hand, the remaining structure
(the null-plane quantum Poincar\'e algebra recently introduced
in \cite{Beyond,Null}) strongly differs from the previous ones:
firstly, it is constructed in a null-plane context where the
Poincar\'e invariance splits into a kinematical and dynamical
part \cite{Dirac} and, secondly, this case is a quantization of
a non-standard (triangular) coboundary Lie bialgebra.

We also recall that the related problem of obtaining 
universal $R$-matrices for  standard  Poincar\'e deformations
has been only solved in (2+1) dimensions \cite{Ita,Itab}.  In
the non-standard case, relevant successes have been recently
obtained: the universal $R$-matrix for the (1+1) case has been
deduced in \cite{Irana,RR,TT} and a (2+1) solution has been
recently given in \cite{Iranb}.

The aim of this letter is twofold: on one hand, to  present a
simplified construction of the universal $R$-matrices for
non-standard quantum $so(2,2)$ and (2+1) Poincar\'e algebras,
which is based on the $sl(2,\R)$ factorized
$R$-matrix introduced in \cite{sl} (sections 2 and 3). 
Whithin this construction a (non-linear) change of basis (whose
origin lies in a $T$-matrix construction
\cite{TT,fg}) is essential, and will lead to rather  new
expressions for all these quantizations. From a physical point
of view, it is important to stress that
$so(2,2)$ is interpreted in a conformal context: i.e.,  its
classical counterpart acts as the group of conformal
transformations on the (1+1) Minkowskian space-time. The second
objective is to get a deeper insight in the (2+1) null-plane
quantization by constructing the associated quantum group and
by exploring its representation theory. All these results are
presented in Section 4.


\sect {Non-standard quantum $so(2,2)$ revisited}

Let us consider the coproduct and the commutation relations of
the non-standard quantum $sl(2,\R)$, denoted by 
$U_z sl(2,\R)=\langle A,A_+,A_-\rangle$, as \cite{sl}:
\bea
&&  \Delta (A_+ ) =1 \otimes A_+  + A_+ \otimes 1,\cr
&& \Delta (A) =1 \otimes A + A\otimes
e^{2z A_+ }, \label{aa} \\
&& \Delta (A_-) = 1 \otimes A_- + A_-\otimes
e^{2 z A_+ }, 
\nonumber
\eea
\be 
  [A,A_+ ]= \frac{e^{2 z A_+} -1  } z,\quad 
[A,A_-]=-2 A_- +z A^2,\quad  [A_+  ,A_- ]= A,
\label{ab} 
\ee 
and the quantum Casimir belonging to the centre  of $U_z
sl(2,\R)$   given by:
\be
{\cal  C}_z=\frac 12 A\, e^{-2zA_+}A+\frac 
{1-e^{-2zA_+}}{2z}A_-
+A_-\frac {1-e^{-2zA_+}}{2z} + e^{-2zA_+}-1 ,
\label{ac}
\ee
These relations are obtained from the ones  given in
\cite{Beyond,Ohn} in terms of
$\{J_3,J_+,J_-\}$, by means of the change of basis \cite{sl}
\bea
&&A_+=J_+,\qquad A= e^{z J_+ }J_3,\cr
&&A_-= e^{z J_+ }J_- - \frac z4 e^{z J_+ } \sinh( z J_+  ).
\label{hh}
\eea

In \cite{sl} it is shown that the universal element
\be
R_z=\exp\{-z A_+\otimes A\}\exp\{z A\otimes A_+\}  
\label{ad}
\ee
is a solution of the quantum Yang--Baxter equation  and also
verifies the property
\be
\sigma \circ \Delta( X) ={ \cal  R}_z \Delta( X) 
{\cal  R}_z^{-1},  \qquad \forall 
\, X \in U_z sl(2,\R),
\label{aac}
\ee
being $\sigma$ the flip operator  $\sigma(a\otimes b)=b\otimes
a$. Hence, $R_z$ is the quantum universal $R$-matrix for $U_z
sl(2,\R)$.

Two comments concerning this universal  $R$-matrix are in
order: firstly,  the significant simplification  obtained (as
far as the commutation rules (\ref{ab}) are concerned) with
respect to the original formulation of this non-standard
$sl(2,\R)$ deformation. Secondly, we recall that the factorized
expression (\ref{ad}) comes from a universal
$T$-matrix formalism \cite{TT,fg}. From this  point of view,
the interest of finding such a kind of factorized expressions
is directly related to the interpretation of the transfer
matrices as quantum monodromies and the obtention of more
manageable algebraic models in quantum field theory (see
\cite{Fr}). 

Let us now consider two copies of the non-standard  quantum
$sl(2,\R)$ algebra, the former with $z$ and the latter with
$-z$ as deformation parameters:
$U_z^{(1)} sl(2,\R)=\langle A^1,A_+^1,A_-^1\rangle$ and 
$U_{-z}^{(2)} sl(2,\R)=\langle A^2,A_+^2,A_-^2\rangle$.  The
generators defined by
\bea
&& \jj=\frac 12(A^1- A^2),\qquad \dd=\frac 12(A^1+A^2),\cr
&& \pu=   A_+^1 + A_+^2 ,\qquad \pd= A_+^1 - A_+^2,\label{ae}\\
&& \cu=   - A_-^1 - A_-^2 ,\qquad \cd= A_-^1 - A_-^2,\nonumber
\eea
give rise to a non-standard quantum deformation of  $so(2,2)$
\cite{Beyond}: 
$$
U_z so(2,2)\simeq U_z^{(1)}  sl(2,\R) \oplus U_{-z}^{(2)} 
 sl(2,\R).
\nonumber
$$
 
At a purely classical level, $SO(2,2)$ can  be regarded in this
basis as the group of conformal transformations of the (1+1)
Minkowskian space-time, where $\jj$ generates the boosts,
$\pu$   the time translations, $\pd$ the space translations,
$\dd$ is a dilation generator and $\cu,\cd$ generate  specific 
conformal transformations. The Hopf algebra structure  of $U_z
so(2,2)$ obtained in this way is   given by the following 
coproduct
$(\Delta)$, counit $(\epsilon)$, antipode $(\gamma)$  and
commutation relations: 
\bea
&&  \Delta (\pu)  =1 \otimes \pu  + \pu \otimes 1,\qquad
\Delta (\pd)  =1 \otimes \pd  + \pd \otimes 1,\cr
 && \Delta (\jj) =1 \otimes \jj + \jj \otimes
e^{z \pd }\cosh z\pu + \dd \otimes e^{z \pd }\sinh z\pu ,\cr
 && \Delta (\dd) =1 \otimes \dd + \dd \otimes
e^{z \pd }\cosh z\pu + \jj \otimes e^{z \pd }\sinh z\pu,
\label{af}\\
 && \Delta (\cu) =1 \otimes \cu + \cu \otimes
e^{z \pd }\cosh z\pu -  \cd \otimes e^{z \pd }\sinh z\pu,\cr
&& \Delta (\cd) =1 \otimes \cd + \cd \otimes
e^{z \pd }\cosh z\pu -  \cu \otimes e^{z \pd }\sinh z\pu,
\nonumber
\eea
\be
\epsilon(X) =0,\qquad \mbox{for $X\in
\{\jj,\pu,\pd,\cu,\cd,\dd\}$},\label{ag}
\ee
\bea
&&\gamma(\pu)=-\pu,\qquad \gamma(\pd)=-\pd,\cr
&&\gamma(\jj)=-\jj  
e^{- z \pd }\cosh z\pu + \dd   e^{- z \pd }\sinh z\pu,\cr
&&\gamma(\dd)=-\dd  
e^{- z \pd }\cosh z\pu + \jj   e^{- z \pd }\sinh z\pu,\label{ah}\\
&&\gamma(\cu)=-\cu  
e^{- z \pd }\cosh z\pu - \cd   e^{- z \pd }\sinh z\pu,\cr
&&\gamma(\cd)=-\cd  
e^{- z \pd }\cosh z\pu - \cu   e^{- z \pd }\sinh z\pu,\nonumber
\eea
\bea
&&[\jj,\pu]=\frac 1z(e^{z\pd} \cosh z\pu -1 ),\qquad 
[\jj,\pd]=\frac 1z e^{z\pd} \sinh z\pu  ,\cr
&&[\jj,\cu]=\cd - z (\jj^2+ \dd^2),\qquad 
[\jj,\cd]=\cu + 2 z \jj\dd  ,\cr
&&[\dd,\pu]=\frac 1z e^{z\pd} \sinh z\pu ,\qquad \quad
[\dd,\pd]=\frac 1z(e^{z\pd} \cosh z\pu -1 )  ,\label{ai}\\
&&[\dd,\cu]=- \cu - 2 z \jj\dd,\qquad\ \quad 
[\dd,\cd]=- \cd + z (\jj^2+ \dd^2) ,\cr
&&[\pu,\cu]=- 2\dd,\quad [\pu,\cd]= 2\jj,\quad 
[\pd,\cu]=-2\jj,\quad
[\pd,\cd]=2\dd,\cr
&&[\jj,\dd]=0,\qquad [\pu,\pd]=0,\qquad [\cu,\cd]=0.\nonumber
\eea
Note that (\ref{af}) presents an interesting feature: both 
generators $\pu$ and $\pd$ are primitive ones.  Therefore, this
conformal approach to $so(2,2)$ leads to a quantum structure
that can be interpreted as an attemp in order to  deform the
(1+1) Minkowskian space and time in a rather symmetrical way.

Two elements of the centre of $U_z so(2,2)$ are  constructed
from the quantum Casimirs (\ref{ac}) of $U_z^{(1)} sl(2,\R)$
and $U_{-z}^{(2)} sl(2,\R)$ as:
\be
{\cal C}^q_1= {\cal C}_z^{(1)} + {\cal C}_{-z}^{(2)},
\qquad
{\cal C}^q_2= {\cal C}_z^{(1)} - {\cal C}_{-z}^{(2)} .
\label{aj}
\ee
After a straightforward computation  we get:
\bea
&&\!\!\!\!\!
\!\!\!\!\!\!{\cal C}^q_1=\jj e^{-z\pd}\cosh (z\pu) \jj + 
\dd e^{-z\pd}\cosh (z\pu) \dd \cr
&&\!\!\!\!\!
\!\!\!\!\!\!\qquad - \jj e^{-z\pd}\sinh (z\pu) \dd
- \dd e^{-z\pd}\sinh (z\pu) \jj\cr
&&\!\!\!\!\!
\!\!\!\!\!\!\qquad +\frac 1{2z}(1-e^{-z\pd}\cosh z\pu)\cd 
+\cd \frac 1{2z}(1-e^{-z\pd}\cosh z\pu)\cr
&&\!\!\!\!\!
\!\!\!\!\!\!\qquad -  e^{-z\pd}\frac {\sinh z\pu}{2z}\cu -
\cu  e^{-z\pd}\frac {\sinh z\pu}{2z}+ 2(e^{-z\pd}\cosh z\pu-1),
\label{ak}\\
&&\!\!\!\!\!
\!\!\!\!\!\! {\cal C}^q_2=\jj e^{-z\pd}\cosh (z\pu) \dd + 
\dd e^{-z\pd}\cosh (z\pu) \jj \cr
&&\!\!\!\!\!
\!\!\!\!\!\!\qquad - \jj e^{-z\pd}\sinh (z\pu) \jj
- \dd e^{-z\pd}\sinh (z\pu) \dd\cr
&&\!\!\!\!\!
\!\!\!\!\!\!\qquad -\frac 1{2z}(1-e^{-z\pd}\cosh z\pu)\cu 
-\cu \frac 1{2z}(1-e^{-z\pd}\cosh z\pu)\cr
&&\!\!\!\!\!
\!\!\!\!\!\!\qquad + e^{-z\pd}\frac {\sinh z\pu}{2z} \cd +
\cd  e^{-z\pd}\frac {\sinh z\pu}{2z} - 2e^{-z\pd}\sinh z\pu .
\label{al}
\eea

Likewise, the universal $R$-matrix for $U_z so(2,2)$ can  be
easily deduced as a product of those corresponding to the two
copies of $U_z  sl(2,\R)$ (\ref{ad}):
\bea
&&\!\!\!\!\!
\!\!\!\!\!\!{\cal R}_z=R_z^{(1)}R_{-z}^{(2)}=
\exp\{-z A^1_+\otimes A^1\}
\exp\{z A^1\otimes
A^1_+\}
\exp\{z A^2_+\otimes A^2\}\exp\{-z A^2\otimes A^2_+\}\cr
&&\quad =  \exp\{-z A^1_+\otimes A^1 + z A^2_+\otimes A^2\}
\exp\{z A^1\otimes A^1_+ -z A^2\otimes A^2_+ \}\cr
&&\quad =  \exp\{-z (\pu \otimes \jj + \pd \otimes \dd)\}
\exp\{z (\jj\otimes \pu + \dd\otimes \pd) \},\nonumber
\eea
which  can be finally written in a complete  ``factorized" form
as:
\be 
{\cal R}_z=\exp\{-z \pu \otimes \jj  \}
\exp\{-z   \pd \otimes \dd \}
\exp\{z   \dd\otimes \pd  \}
\exp\{z  \jj\otimes \pu  \}.
\label{am}
\ee
The first order in $z$ gives the classical $r$-matrix
\be
r=z(\jj\wedge \pu + \dd\wedge\pd),
\label{an}
\ee
which, as expected, is a solution of the classical Yang--Baxter 
equation.


\sect {(2+1) null-plane quantum Poincar\'e algebra}

\subsect {Null-plane classical  Poincar\'e algebra}

We  briefly describe the classical structure of the (2+1)
Poincar\'e algebra ${\cal P}(2+1)$ in relation with  the
null-plane evolution scheme \cite{LS}, in which the initial
state of a quantum relativistic system can be defined on a
light-like plane $\Pi_n^\tau$ defined by $n\cdot x=\tau$, where
$n$ is a light-like vector and $\tau$ a real constant.  In
particular, if 
$n=(\frac 12,0,\frac 12)$ and the coordinates
\be
x^- = n\cdot x =\frac 12(x^0-x^2),\qquad
x^+=x^0+x^2, \label{ba} 
\ee
are considered,   
a point $x\in \Pi_n^\tau$ 
will be labelled by $(x^+, x^1)$ while the
remaining one ($x^-$) plays the role of a time parameter $\tau$. 
A basis $\{P_+,P_1,P_-,E_1,F_1,K_2\}$ of the (2+1)  Poincar\'e
algebra  consistent with these coordinates 
is provided by the   generators $P_+$, $P_-$, $E_1$ and
$F_1$ which are defined in the terms of the usual  kinematical
ones
$\{P_0,P_1,P_2,K_1,K_2,J\}$ by:
\be
{P_+}=\frac 12(P_0+P_2),\quad {P_-}=P_0-P_2,\quad 
{E_1}=\frac 12(K_1+J),\quad {F_1}=K_1-J .
\label{bb}
\ee
This ``null-plane" basis has the following  non-vanishing
commutation rules:
\bea
&&[K_2,P_\pm]=\pm P_\pm,\quad 
 [K_2,E_1]=E_1   ,\quad 
[K_2,F_1]=- F_1  ,\cr
&&[E_1,P_1]= P_+ ,\quad 
[F_1,P_1]=  P_-  ,\quad
 [E_1,F_1]=K_2,\label{bc}\\
&&[P_+,F_1]=  -  P_1,\quad 
[P_-,E_1]=  -  P_1  .
\nonumber
\eea
The operators $\{P_+, P_1,  E_1,K_2\}$ are  the infinitesimal
generators of the stability group $S_+$ of the null-plane
$\Pi_n^0$ $(\tau=0)$. The remaining generators have a dynamical
significance: the hamiltonian $P_-$ translates $\Pi_n^0 $ into
$\Pi_n^\tau$ and $F_1$ rotates it around the surface of the
light-cone $x^1=0$.
 
Finally, let us recall that the centre of  ${\cal P}(2+1)$ is
generated by the square of the mass operator $M^2$ and the
intrinsic angular momentum $L$, which read:
\bea
&& M^2=2P_-P_+-P_1^2,\label{bd}\\
&&L=K_2P_1+E_1P_- -F_1P_+ .
\label{be}
\eea

\subsect {Null-plane quantum  Poincar\'e algebra}

This null-plane Poincar\'e algebra is naturally  linked by a
contraction procedure to $so(2,2)$ when the latter is written
in a basis of the kind (\ref{ae}).
The explicit form of that contraction mapping is as follows:
\bea
&& {P_+}=\varepsilon \, \frac 1{\sqrt{2}} \pd,\quad 
{P_1}=\varepsilon\, \jj,\quad 
{P_-}=- \varepsilon\, \frac 1{\sqrt{2}} \cd,\cr
&&{E_1}=-\frac 1{\sqrt{2}} \pu,\quad 
{F_1}= \frac 1{\sqrt{2}} \cu,\quad {K_2}=D.
\label{cb}
\eea
In the quantum case, the deformation parameter  has also to be
transformed as  $w= \frac 1{\varepsilon\sqrt{2}}
z$. Therefore, 
by applying (\ref{cb}) in the Hopf algebra of $U_z so(2,2)$
(\ref{af}--\ref{ai}) and then by making the limit 
$\varepsilon\to 0$ we get the resulting Hopf structure for the
quantum (2+1) null plane Poincar\'e algebra
$U_w{\cal P}(2+1)$:
\bea
&&  \Delta (P_+)  =1 \otimes P_+  + P_+ \otimes 1,\qquad
\Delta (E_1)  =1 \otimes E_1  + E_1 \otimes 1,\cr
&& \Delta (P_-) =1 \otimes P_- + P_- \otimes e^{2 w P_+ },
\qquad \Delta (P_1) =1 \otimes P_1
+ P_1 \otimes e^{2 w P_+ },\cr
&& \Delta (F_1) =1 \otimes F_1
+ F_1 \otimes e^{2 w P_+ } - 2 w P_-\otimes e^{2 w P_+ }E_1,
\label{cc}\\
&& \Delta (K_2) =1 \otimes K_2
+ K_2 \otimes e^{2 w P_+ } - 2 w P_1\otimes e^{2 w P_+ }E_1,
\nonumber
\eea
\be
\epsilon(X) =0,\qquad \mbox{for $X\in
\{{K_2},{P_+},{P_-},{P_1},{E_1},{F_1}\}$},
\ee
\bea
&&\gamma(P_+)=-P_+,\qquad \gamma(E_1)=-E_1,\cr
&& \gamma(P_-)=-P_-e^{-2wP_+}, \quad 
\gamma(P_1)=-P_1e^{-2wP_+} , \\
&&\gamma(F_1)=-F_1e^{-2wP_+}-2w P_-e^{-2wP_+}E_1,\cr
&&\gamma(K_2)=-K_2e^{-2wP_+}-2w P_1e^{-2wP_+}E_1,\nonumber
\eea
\bea
&&[K_2,P_+]=\frac 1{2w}(e^{2wP_+}  -1 ),\qquad 
[K_2,P_-]=- P_- - wP_1^2 ,\cr
&&[K_2,E_1]=E_1 e^{2wP_+}  ,\qquad 
[K_2,F_1]=- F_1 - 2wP_1 K_2 ,\label{cd}\\
&&[E_1,P_1]=\frac 1{2w}(e^{2wP_+}  -1 ),\qquad 
[F_1,P_1]=  P_- + wP_1^2 ,\cr
&&[E_1,F_1]=K_2,\qquad 
[P_+,F_1]=  -  P_1,\qquad 
[P_-,E_1]=  -  P_1  ,
\nonumber
\eea
where the remaining commutators are zero.
Note that the generators of the null-plane  stability group
close a Hopf subalgebra   $U_wS_+$. As a byproduct of the
original change of basis within
$sl(2,\R)$, these commutation rules are simpler  than the ones
given in
\cite{Beyond}.

The quantum Casimirs belonging to the  centre of $U_w{\cal
P}(2+1)$ are deduced from (\ref{ak}) and (\ref{al}) by means of
the limits:
\be
M_q^2=\lim_{\varepsilon\to 0}(-\varepsilon^2 {\cal C}_1^q),
\qquad
L_q=\frac 12\lim_{\varepsilon\to 0}(\varepsilon  {\cal C}_2^q),
\ee
explicitly,
\bea
&&M_q^2=P_-\frac{1-e^{-2wP_+}}{w} -P_1^2e^{-2wP_+},\label{ce}\\
&&L_q=K_2P_1e^{-2wP_+}+E_1(P_-+wP_1^2)e^{-2wP_+}-F_1
\frac{1-e^{-2wP_+}}{2w}
.\label{cf}
\eea
The universal $R$-matrix $U_w{\cal P}(2+1)$  is also directly
obtained from (\ref{am}) and reads: 
\be 
{\cal R}_w=\exp\{2 w E_1 \otimes P_1  \}
\exp\{-2w   P_+ \otimes K_2 \}
\exp\{2w   K_2\otimes P_+  \}
\exp\{-2w  P_1\otimes E_1  \}.
\label{cg}
\ee

A   differential representation of $U_wS_+$  with coordinates
$(p_+,p_1)$ can be given as follows:
\be 
 P_+=p_+,\quad P_1=p_1,\quad K_2=
\frac{e^{2wp_+}-1}{2w}\partial_+,\quad 
  E_1=\frac{e^{2wp_+}-1}{2w}\partial_1,
\label{ch}
\ee 
where $\partial_+=\frac{\partial}{\partial p_+}$ and 
$\partial_1=\frac{\partial}{\partial p_1}$.  With the aid of
the quantum Casimirs a spin-zero differential representation
$(L_q=0)$ for the two remaining generators of $U_w{\cal
P}(2+1)$ can be deduced:
\be 
 P_-=\frac{w(m_q^2+p_1^2e^{-2wp_+})}{1-e^{-2wp_+}},\quad
 F_1=p_1\partial_++
\frac{w(m_q^2+p_1^2)}{1-e^{-2wp_+}}\partial_1,
\label{ci}
\ee 
where $m_q^2$ is the eigenvalue of the  $q$-Casimir (\ref{ce}).
Similarly to the classical case, we can take the coordinate
$x^-$ as an evolution parameter ($\tau$) and thus we can
consider a wave function $\psi(p_+,p_1,\tau)$ whose evolution
is determined by   the $q$-Scr\"odinger equation provided by
the Hamiltonian $P_-$: $i\partial_\tau\psi=P_-\psi$. In terms
of the representation (\ref{ci}) we get: \be
i\partial_\tau\psi(p_+,p_1,\tau)=
\frac{w(m_q^2+p_1^2e^{-2wp_+})}{1-e^{-2wp_+}}
\psi(p_+,p_1,\tau),
\label{cj}
\ee
which is different from the one  given in \cite{Null} for the
(3+1) case. This fact can be more clearly appreciated by
writing the power series expansion in $w$ of  $P_-$:
\bea
&&P_-=\frac{w(m_q^2+p_1^2e^{-2wp_+})}{1-e^{-2wp_+}}
=\frac{w(m_q^2e^{wp_+}+p_1^2e^{-wp_+})}{2\sinh wp_+}\cr
&&=\frac{  m_q^2+p_1^2 }{2p_+}+w\frac{ m_q^2-p_1^2 }{2}
+w^2\frac{p_+(  m_q^2+p_1^2) }{6}+o(w^3).
\label{ck}
\eea
The zero-term in $w$ can be  identified with a kinetic term  of
the null-plane bound state equation in quantum chromodynamics
\cite{KSop,HHLS} while   all remaining terms in $w$ constitute
a dynamical part, now including a first order term in $w$ (that
is absent in \cite{Null}). Therefore, this deformation of the
null-plane symmetry has some intrinsic dynamical content whose
explicit description depends on the way in which the
deformation is constructed.


\sect {(2+1) Null-plane quantum  Poincar\'e group}

The Lie bialgebra underlying the quantum Hopf  algebra of
$U_w{\cal P}(2+1)$ is generated by the non-standard classical 
$r$-matrix (first order in $w$ of (\ref{cg})):
\be
r=2(K_2\wedge P_+ + E_1\wedge P_1) ,
\label{bbe}
\ee
 which provides   the cocommutators 
$\delta(X)=[1\otimes X+X\otimes 1,r]$:
\bea 
&& \delta(P_+)=0,\qquad \delta(E_1)=0, \cr 
&&\delta(P_1)=2  P_1\wedge P_+ ,\qquad 
\delta(P_-)=2  P_-\wedge P_+
,\label{bbf}\\ 
&&\delta(F_1)=2 (F_1\wedge P_+ +E_1\wedge P_-),\cr 
&&\delta(K_2)=2 (K_2\wedge P_+ +E_1\wedge P_1 ).\nonumber  
\eea
They are related to   the first order term in  the deformation
parameter of the coproduct (\ref{cc})  by means of
$\delta=\Delta_{(1)} -\sigma\circ \Delta_{(1)}$.

The $r$-matrix (\ref{bbe}) also allows to deduce  the
associated Poisson structure to the Poincar\'e algebra. Let the
four-dimensional matrix representation of ${\cal P}(2+1)$  
given by: 
\bea
 D(P_+)=\left(\begin{array}{cccc}
 0 & 0 & 0 & 0 \\
\frac 12 & 0 & 0 & 0 \\
0 & 0 & 0 & 0 \\
\frac 12 & 0 & 0 & 0 \end{array}\right)\  
D(P_-)=\left(\begin{array}{cccc}
 0 & 0 & 0 & 0 \\
1 & 0 & 0 & 0 \\
0 & 0 & 0 & 0 \\
-1 & 0 & 0 & 0 \end{array}\right)&&\!\!\!\!\!\! 
D(P_1)=\left(\begin{array}{cccc}
 0 & 0 & 0 & 0 \\
0 & 0 & 0 & 0 \\
1 & 0 & 0 & 0 \\
0 & 0 & 0 & 0 \end{array}\right) \cr
 D(E_1)=\left(\begin{array}{cccc}
 0 & 0 & 0 & 0 \\
0 & 0 & \frac 12 & 0 \\
0 & \frac 12 & 0 & -\frac 12 \\
0 & 0 & \frac 12 & 0 \end{array}\right)\ 
D(F_1)=\left(\begin{array}{cccc}
 0 & 0 & 0 & 0 \\
0 & 0 & 1 & 0 \\
0 & 1 & 0 & 1 \\
0 & 0 & -1 & 0 \end{array}\right)&&\!\!\!\!\!\! 
D(K_2)=\left(\begin{array}{cccc}
 0 & 0 & 0 & 0 \\
0 & 0 & 0 & 1 \\
0 & 0 & 0 & 0 \\
0 & 1 & 0 & 0 \end{array}\right)\cr
&& \label{bf} 
\eea
Then a $4\times 4$ representation of the element 
$g=e^{a^+P_+}e^{a^-P_-}e^{a^1P_1}e^{e^1E_1}
e^{f^1F_1}e^{k^2K_2}$ belonging to
the   (2+1) Poincar\'e group is
\be
D(g)=\left(\begin{array}{cccc}
 1 & 0 & 0 & 0 \\
\frac {a^+}2+a^- & \Lambda_0^0 & \Lambda_1^0 & \Lambda_2^0 \\
a^1 & \Lambda_0^1 & \Lambda_1^1 & \Lambda_2^1 \\
\frac {a^+}2-a^- & \Lambda_0^2 & \Lambda_1^2 & \Lambda_2^2
 \end{array}\right),
\label{bg} 
\ee
where the $\Lambda_\nu^\mu$ are the matrix  elements of the
Lorentz subgroup
 (whose generators are $E_1$, $F_1$ and $K_2$) satisfying the
pseudo-orthogonality condition:  
\be
\Lambda_\nu^\mu\Lambda_\sigma^\rho\eta^{\nu\sigma}=
\eta^{\mu\rho},
\quad (\eta^{\mu\rho})={\mbox{diag}}\,(1,-1,-1).
\label{bh}
\ee
The Poisson brackets
of the coordinate functions on the Poincar\'e group  are
obtained by calculating the Poisson bivector
\be
\{D(g)\dot\otimes D(g)\}=[r,D(g)\dot\otimes D(g)],
\label{bi}
\ee
writing the $r$-matrix (\ref{bbe}) in terms of the matrix
representation (\ref{bf}).
The final result can be summarized as follows:
\bea
&&\{ { a}^+,{ a}^1\}=-2w\,{ a}^1,\quad 
\{ { a}^+,{ a}^-\}=-2w\,{ a}^-
,\quad  \{ { a}^1,{ a}^-\}=0,\label{bj}\cr 
&&\{ {\Lambda}_\nu^\mu,{\Lambda}_\rho^\sigma\}=0,\qquad
\nu,\mu,\rho,\sigma=0,1,2;\cr
&&\{ {\Lambda}_\nu^\mu,{ a}^+\}=-2\delta_{\mu 0}
{\Lambda}_\nu^2
-2\delta_{\mu 2}{\Lambda}_\nu^0-(\mu-1)(\nu-1)+(\Lambda_0^\mu
+\Lambda_2^\mu)(\Lambda_\nu^0+\Lambda_\nu^2),\cr
&&\{ {\Lambda}_\nu^\mu,{ a}^1\}= \delta_{\mu 1}(1-\nu-
\Lambda_\nu^0+\Lambda_\nu^1+\Lambda_\nu^2)+
\Lambda_\nu^1(\Lambda_0^\mu+\Lambda_2^\mu-1),\cr
&&\{ {\Lambda}_\nu^\mu,{ a}^-\}=\frac 12 (\mu-1)^2(\nu-1)+
\frac 12  (\Lambda_0^\mu
+\Lambda_2^\mu)(\Lambda_\nu^0-\Lambda_\nu^2).
\label{bbf}
\eea
It is worth comparing these expressions with the  results
related to the classical $r$-matrix of the $\kappa$-Poincar\'e
algebra \cite{Zak,Mas,Lukdos}.

The classical matrix  representation (\ref{bf})  is also valid
for $U_w{\cal P}(2+1)$ since $D(P_+)^2$ vanishes. 
This fact can be used to get an explicit expression for 
${\cal R}_w$: 
\be 
  D({\cal R}_w)=I\otimes I +2w(D(K_2)\wedge D(P_+) 
+ D(E_1)\wedge D(P_1)),
\label{di}
\ee
where $I$ is the
$4\times 4$ identity matrix. The fulfillment of property
(\ref{aac}) allows to apply the FRT method \cite{FRT}:
\be
RT_1T_2=T_2T_1R,
\label{dj}
\ee
where $R$ is (\ref{di}), $T_1=T\otimes I$,  $T_2=I\otimes T$,
being $T$ the group element (\ref{bg}) but now with
non-commutative entries:
${\hat\Lambda}_\nu^\mu$ and ${\hat a}^i$. The commutation
relations of the quantum Poincar\'e group read
\bea
&&[ {\hat a}^+,{\hat  a}^1]=-2w\,{\hat  a}^1,\quad 
[ {\hat  a}^+,{\hat  a}^-]=-2w\,{\hat  a}^-
,\quad  [ {\hat  a}^1,{\hat  a}^-]=0,\label{dk}\cr 
&&[ {\hat \Lambda}_\nu^\mu,{\hat \Lambda}_\rho^\sigma]=0,\qquad
\nu,\mu,\rho,\sigma=0,1,2;\cr
&&[ {\hat \Lambda}_\nu^\mu,{\hat  a}^+]= -2\delta_{\mu 0}{\hat
\Lambda}_\nu^2 -2\delta_{\mu 2}{\hat
\Lambda}_\nu^0-(\mu-1)(\nu-1)+({\hat  \Lambda}_0^\mu +{\hat
\Lambda}_2^\mu)({\hat \Lambda}_\nu^0+{\hat \Lambda}_\nu^2),\cr
&&[ {\hat \Lambda}_\nu^\mu,{\hat  a}^1]= \delta_{\mu 1}(1-\nu-
{\hat \Lambda}_\nu^0+{\hat \Lambda}_\nu^1+{\hat
\Lambda}_\nu^2)+ {\hat \Lambda}_\nu^1({\hat
\Lambda}_0^\mu+{\hat \Lambda}_2^\mu-1),\cr &&[ {\hat
\Lambda}_\nu^\mu,{\hat  a}^-]=\frac 12 (\mu-1)^2(\nu-1)+
\frac 12  ({\hat \Lambda}_0^\mu
+{\hat \Lambda}_2^\mu)({\hat \Lambda}_\nu^0- {\hat
\Lambda}_\nu^2),
\label{ddf}
\eea
with the additional relations:
\be
{\hat \Lambda}_\nu^\mu{\hat \Lambda}_\sigma^\rho
\eta^{\nu\sigma}=\eta^{\mu\rho},
\quad (\eta^{\mu\rho})={\mbox{diag}}\,(1,-1,-1).
\label{dl}
\ee
As it happened   with   the $\kappa$-Poincar\'e  group
\cite{Mas,Lukdos} these commutation relations  are also a Weyl
quantization $[\, ,\,]\to w^{-1}[\, ,\,]$  of the Poisson
brackets of the coordinate functions on the Poincar\'e group
(\ref{bbf}), and moreover, all  the Lorentz coordinates
${\hat \Lambda}_\nu^\mu$ commute among themselves so that 
there is no ordering
ambiguity.  The associated coproduct, counit and
antipode can be   deduced from relations
$\Delta(T)=T\dot\otimes T$, $\epsilon(T)=I$  and
$\gamma(T)=T^{-1}$, respectively. In particular, the coproduct
is:
\bea
&&\Delta({\hat a}^+)={\hat a}^+\otimes 1 +
\frac 12 ({\hat \Lambda}_0^0+{\hat \Lambda}_0^2
+{\hat \Lambda}_2^0+{\hat \Lambda}_2^2)\otimes {\hat a}^+\cr
&&\qquad + ({\hat \Lambda}_1^0+
{\hat \Lambda}_1^2)\otimes  {\hat a}^1
+({\hat \Lambda}_0^0+{\hat \Lambda}_0^2
-{\hat \Lambda}_2^0-{\hat \Lambda}_2^2)\otimes {\hat a}^-,\cr
&&\Delta({\hat a}^1)={\hat a}^1\otimes 1 +
\frac 12 ({\hat \Lambda}_0^1 +  {\hat \Lambda}_2^1 )
\otimes {\hat a}^+
+{\hat \Lambda}_1^1\otimes {\hat a}^1 +
 ({\hat \Lambda}_0^1-{\hat \Lambda}_2^1) \otimes  {\hat a}^- 
,\label{dm}\\ 
&&\Delta({\hat a}^-)={\hat a}^-\otimes 1 +
\frac 14
({\hat \Lambda}_0^0-{\hat \Lambda}_0^2
+{\hat \Lambda}_2^0-{\hat \Lambda}_2^2)\otimes  {\hat a}^+ \cr
&&\qquad + 
\frac 12 ({\hat \Lambda}_1^0-{\hat \Lambda}_1^2)\otimes {\hat
a}^1  +\frac 12
({\hat \Lambda}_0^0-{\hat \Lambda}_0^2
-{\hat \Lambda}_2^0+{\hat \Lambda}_2^2)\otimes {\hat a}^- ,\cr
&&\Delta({\hat \Lambda}_\nu^\mu)=
{\hat \Lambda}_\sigma^\mu\otimes
{\hat \Lambda}_\nu^\sigma .
\nonumber
\eea
The quantum (2+1) Poincar\'e plane of coordinates $({\hat
x}^+,{\hat x}^1, {\hat x}^-)$  characterized by
\be
[{\hat x}^+,{\hat x}^1]=-2w\,{\hat x}^1,
\quad [{\hat x}^+,{\hat x}^-]=-2w\,{\hat
x}^- ,\quad  [{\hat x}^1,{\hat x}^-]=0,
\label{dn}
\ee
is easily derived from  the first three commutators of
(\ref{dk}). Note that it includes, as a particular case, the
quantum (1+1) Poincar\'e plane $[{\hat x}^+,{\hat
x}^-]=-2w\,{\hat x}^-$ \cite{TT}. The coordinates 
$({\hat x}^+,{\hat x}^1)$ could be interpreted as the
parameters of a quantum light-like plane while the remaining
one ${\hat x}^-$ would be a quantum time.



\bigskip

\noindent
{\large{{\bf Acknowledgements}}}

\bigskip

 This work has been
partially supported by DGICYT (Project PB94--1115) from the
Ministerio de Educaci\'on y Ciencia de Espa\~na.

\bigskip


\end{document}